\documentclass[12pt,a4paper]{amsart}

\usepackage{amssymb}
\makeatletter
\renewcommand{\@begintheorem}[2]{
\rm \trivlist \item [\hskip \labelsep {\bf #2\ \ #1.}]
                                }
\makeatother

\makeatletter

\DeclareFontFamily{U}{cyr}{}
\DeclareFontShape{U}{cyr}{m}{n}{
  <5> wncyr5 <6> wncyr6 <7> wncyr7 <8> wncyr8 <9> wncyr9 <10->
wncyr10}{}
\DeclareMathAlphabet{\mathcyr}{U}{cyr}{m}{n}

\input cyracc.def
\input epsf

\newcommand{\ts}{\vspace{\baselineskip}\noindent{\bf Proof.}$\;\;$}
\newcommand{\ZZ}{{\bf Z}}
\newcommand{\QQ}{{\bf Q}}
\newcommand{\RR}{{\bf R}}
\newcommand{\CC}{{\bf C}}

\newcommand{\HH}{{\bf H}}
\newcommand{\PP}{{\bf P}}

\newcommand{\bA}{{\mathbb A}}

\newcommand{\cA}{{\mathcal A}}

\newcommand{\cG}{{\mathcal G}}

\newcommand{\ccH}{{\mathcal H}}
\newcommand{\cL}{{\mathcal L}}
\newcommand{\cO}{{\mathcal O}}

\newcommand{\im}{\mbox{im}}

\newcommand{\rmd}{\mbox{d}}

\newcommand{\bes}{\begin{equation*}}
\newcommand{\ees}{\end{equation*}}

\title{Abelian surfaces with an automorphism and quaternionic multiplication}
\author{Matteo A. Bonfanti}

\author{Bert van Geemen}
\address{Dipartimento di Matematica, Universit\`a di Milano,
Via Saldini 50, 20133 Milano, Italia}
\email{matteoalfonso.bonfanti@unimi.it}
\email{lambertus.vangeemen@unimi.it}
\begin{document}

\begin{abstract}
We construct one dimensional families of Abelian surfaces with quaternionic multiplication
which also have an automorphism of order three or four. Using Barth's description of the moduli space of $(2,4)$-polarized Abelian surfaces, we find the Shimura curve parametrizing these Abelian surfaces in a specific case. 
We explicitly relate these surfaces to the Jacobians of genus two curves studied by Hashimoto and Murabayashi. 
We also describe a (Humbert) surface in Barth's moduli space which parametrizes Abelian surfaces with real multiplication by $\ZZ[\sqrt{2}]$.
\end{abstract}

\maketitle

\subsection*{Introduction} 
The Abelian surfaces, with a polarization of a fixed type, 
whose endomorphism ring is an order in a quaternion algebra are parametrized by a curve, 
called a Shimura curve, in the moduli space of polarized Abelian surfaces. 
There have been several attempts to find concrete examples of such Shimura curves and of the family of Abelian surfaces over this curve. 
In \cite{HM}, Hashimoto and Murabayashi find two 
Shimura curves as the intersection, 
in the moduli space of principally polarized Abelian surfaces, 
of two Humbert surfaces. 
Such Humbert surfaces are now known `explicitly' also in many other cases (see \cite{BW})
and this might allow one to find explicit models of other Shimura curves. 
Another approach was taken by Elkies in \cite{Elkies} who characterizes elliptic fibrations on the Kummer surfaces of such  Abelian surfaces.

In this paper we consider the rather special case that one of the Abelian surfaces in the family is the selfproduct of an elliptic curve. Moreover we assume
this elliptic curve to have an automorphism (fixing the origin) of order three or four. 
It is then easy to show that, for a fixed product polarization of type $(1,d)$, 
the deformations of the selfproduct with the automorphism are parametrized by a Shimura curve. 
In fact, an Abelian surface with such an automorphism must have a N\'eron Severi group of rank at least three and we show that this implies that the endomorphism algebra 
of such a surface is in general a quaternion algebra. 
One can then work out for which $d$ the quaternion algebra is actually a skew field 
(rather than a matrix algebra). 
The cases for $d\leq 20$ are listed in Section \ref{cordisc}.

The remainder of this paper is devoted to the case of an automorphism of order three and a polarization of type $(1,2)$. In that case the general endomorphism ring is a maximal order $\cO_6$ of the quaternion algebra of discriminant $6$. 
Barth, in \cite{Barth}, provides a description of a moduli space $M_{2,4}$, 
embedded in $\PP^5$, of $(2,4)$-polarized Abelian surfaces with a level structure. 
Since the polarized Abelian surfaces we consider have an automorphism of order three, 
the corresponding points in $M_{2,4}$ are fixed by an automorphism of order three.
This allows us to explicitly identify the Shimura curve in $M_{2,4}$ which 
parametrizes the Abelian surfaces 
with quaternionic multiplication by the maximal order $\cO_6$ 
in the quaternion algebra with discriminant $6$. 
It is embedded as a line, which we denote by $\PP^1_{QM}$, in $M_{2,4}\subset\PP^5$. 
The symmetric group $S_4$ acts on this line by changing the level structures.

According to Rotger \cite{ Rotger}, an Abelian surface with endomorphism ring $\cO_6$
has a unique principal polarization, 
which is in general defined by a genus two curve in that surface.
We show explicitly how to find such genus two curves, or rather their images in the Kummer surface embedded in $\PP^5$ with a $(2,4)$-polarization. 
These curves were already considered by Hashimoto and Murabayashi in \cite{HM}. 
We give the explicit relation between the two descriptions in Proposition \ref{isoc}. 
As a byproduct, we find a (Humbert) surface 
in $M_{2,4}$ that parametrizes Abelian surfaces with $\ZZ[\sqrt{2}]$ in the endomorphism ring.
 
In a series of papers (\cite{GPe}, $\ldots$, \cite{GP2}), Gross and Popescu studied, 
both in general and for several small $d$ in particular, explicit maps from 
moduli spaces of $(1,d)$-polarized Abelian surfaces to projective spaces. 
The methods we used to find the Shimura curve in $M_{2,4}$ can, in principle,
be extended also to these cases.

\subsection*{Acknowledgements}
We are indebted to Igor Dolgachev, Alice Garbagnati and Shuang Yan for stimulating
discussions.

\

\section{Polarized Abelian surfaces with automorphisms}\label{one}

\subsection{Abelian surfaces with a $(1,d)$-polarization}\label{aspol}
We recall the basic results on moduli spaces of Abelian surfaces with a $(1,d)$-polarization, following \cite[Chapter 1]{HKW}. Such an Abelian surface $A$ is isomorphic to $\CC^2/\Lambda$, where  the lattice $\Lambda$ can be obtained as the image of 
$\ZZ^4$ under the map given by the period matrix $\Omega$, where we consider all vectors as row vectors:
$$
A\,\cong\,\CC^2/ \Lambda,\qquad
 \Lambda\,=\,\ZZ^4\Omega,\quad
\Omega\,:\,\ZZ^4\,\longrightarrow\,\CC^2,\qquad
x\,\longmapsto\, x\Omega\,=\,x\!\begin{pmatrix}\tau\\ \Delta_d\end{pmatrix}
\,=\,x\!\begin{pmatrix}
\tau_{11}&\tau_{12}\\\tau_{21} &\tau_{22}\\ 
1&0\\0&d\end{pmatrix}~,
$$
where $\tau$ is a symmetric complex $2\times 2$ matrix with positive definite imaginary part, so $\tau\in\HH_2$, the Siegel space of degree two, and $\Delta_d$ is a diagonal matrix with entries $1,d$.
The polarization on $A$ is defined by the Chern class of an ample line bundle in 
$H^2(A,\ZZ)\cong \wedge^2H^1(A,\ZZ)=\wedge^2Hom(\Lambda,\ZZ)$, 
that is, by an alternating map $E_{d}:\Lambda\times\Lambda\rightarrow\ZZ$ 
which is the one defined by the alternating matrix with the same name (so $E_{d}(x,y)=xE_{d}{}^t\!y)$:
$$
E_{d}\,:=\,\begin{pmatrix}0&\Delta_d\\-\Delta_d&0\end{pmatrix}\,=\,
\begin{pmatrix}0&0&1&0\\
0&0&0&d\\-1&0&0&0\\0&-d&0&0\end{pmatrix}~.
$$

\subsection{Products of elliptic curves}\label{prelc}
The selfproduct of an elliptic curve with an automorphism 
of order three and four respectively provides, for any integer $d>0$, a $(1,d)$-polarized Abelian surface with an automorphism of the same order whose eigenvalue on $H^{2,0}$ is equal to one.

To see this,
let $\zeta_j:=e^{2\pi i/j}$ be a primitive $j$-th root of unity. For $j=3,4$, 
let $E_j$ be the elliptic curve with an automorphism $f_j\in \mbox{End}(E_j)$ of order $j$:
$$
E_j\,:=\,\CC\,/\,\ZZ+\ZZ \zeta_j,\qquad 
f_j:E_j\,\longrightarrow\,E_j,\quad z\,\longmapsto\,\zeta_jz~.
$$
Then the Abelian surface $A_j:=E_j^2$ has the automorphism 
$$
\phi_j\,:=\,f_j\times f_j^{-1}:\;A_j:=\,E_j\times E_j\,\longrightarrow\,A_j~.
$$

As $f_j^*$ acts as multiplication by $\zeta_j$ on 
$H^{1,0}(E_j)=\CC\rmd z$, the eigenvalues of $\phi_j^*$ on 
$H^{1,0}(A_j)$ are $\zeta_j,\zeta_j^{-1}$. Thus $\phi_j^*$ acts as the identity on $H^{2,0}(A_j)=\wedge^2H^{1,0}(A_j)$.

The principal polarization on $E_j$ is fixed by $f_j$, so the product of 
this polarization on the first factor with $d$-times the principal polarization on
the second factor is a $(1,d)$-polarization on $A_j$ which is invariant under $\phi_j$.

The lattice $\Lambda_j\subset \CC^2$ defining $A_j$ is given 
by the image of the period matrix $\Omega_j$:
$$
A_j\,\cong\,\CC^2\,/\,\Lambda_j,\qquad \Lambda_j\,=\,\ZZ^4\Omega_j,
\qquad
\Omega_j\,:=\,\begin{pmatrix}\zeta_j&0\\ 0&d\zeta_j\\ 1&0\\0&d\end{pmatrix}~.
$$
The automorphism $\phi_j$ determines, and is determined by, the matrices $\rho_r(\phi_j)$ and $\rho_a(\phi_j)$ which give the action of $\phi_j$ on $\Lambda_j$ and $\CC^2$ respectively. Here we have
$$
\rho_r(\phi_j)\Omega_j\,=\,\Omega_j\rho_a(\phi_j),\qquad
\rho_r(\phi_j)\,\,=\,M_j,\quad
\rho_a(\phi_j)\,=\,\begin{pmatrix} \zeta_j&0\\0&\zeta_j^{-1}\end{pmatrix}~,
$$
where the matrix $M_j$ is given by:
$$
M_3\,:=\,\begin{pmatrix}-1&0&-1&0\\
0&0&0&1\\1&0&0&0\\0&-1&0&-1\end{pmatrix},\qquad
M_4\,:=\,\begin{pmatrix}0&0&-1&0\\
0&0&0&1\\1&0&0&0\\0&-1&0&0\end{pmatrix}~.
$$
The $(1,d)$-polarization is defined by the alternating matrix $E_{d}$, 
from Section \ref{aspol} and is indeed preserved
by $\phi_j$ (so $\phi_j^*E_{d}=E_{d}$) since
$
M_jE_{d}\,{}^t\!M_j\,=\,E_{d}
$.

\subsection{Deformations of $(A_j,E_{1,d},\phi_j)$}\label{defpr}
For a matrix $M\in M_4(\RR)$ such that $ME_d{}^t\!M=E_d$ we define
$$
M\ast_d\tau\,:=\,
(A\tau+B\Delta_d)(C\tau+D\Delta_d)^{-1}\Delta_d,\qquad
\mbox{where}\quad
M\,=\,\begin{pmatrix}A&B\\C&D\end{pmatrix}~.
$$
The  fixed point set of $M_j$ for the $\ast_d$-action on $\HH_2$ is denoted by
$$
\HH_{j,d}\,:=\,\left\{\,\tau\,\in\,\HH_2\,: \; M_j\ast_d\tau\,=\,\tau\,\right\}~.
$$

The following proposition shows that the $(1,d)$-polarized Abelian surfaces
which are deformations of $(A_j,\phi_j)$
form a one parameter family which is parametrized by $\HH_{j,d}$.
We will see in Theorem \ref{enddef} and Corollary \ref{cordisc} that for certain combinations of $j$ and $d$ the general surface in this family is simple and has quaternionic multiplication.

\subsection{Proposition}\label{defm}
The $(1,d)$-polarized Abelian surface  
$(A_{\tau,d}=\CC^2/(\ZZ^4\Omega_\tau),E_{d})$, with $\tau\in\HH_2$, admits an automorphism $\phi_j$ induced by $M_j$ if and only if $\tau\in\HH_{j,d}$.
Moreover, $\HH_{j,d}$ is biholomorphic to  $\HH_1$, the Siegel space of degree one.

\ts
The Abelian surface $A_{\tau,d}=\CC^2/(\ZZ^4\Omega_\tau)$ admits an automorphism induced by $M_j$ if there is a $2\times 2$ complex matrix $N_\tau$ such that 
$$
M_j \Omega_\tau\,=\,\Omega_\tau N_\tau,\qquad
\Omega_\tau\,:=\,\begin{pmatrix}\tau\\ \Delta_{d}\end{pmatrix}~.
$$
Writing $M_j$ as a block matrix with rows $A,B$ and $C,D$, the equation 
$M_j \Omega_\tau\,=\,\Omega_\tau N_\tau$ 
is equivalent to the two equations
$$
A\tau+B\Delta_d\,=\,\tau N_\Omega,\qquad
C\tau+D\Delta_d\,=\,\Delta_d N_\tau~,
$$
hence $N_\tau=\Delta_d^{-1}(C\tau+D\Delta_d)$ and substituting this in the first equation we get:
$$
(A\tau+B\Delta_d)(C\tau+D\Delta_d)^{-1}\Delta_d\,=\,\tau\qquad
\mbox{hence}\quad M_j\ast_d\tau\,=\,\tau~.
$$
Conversely, if 
$ M_j\ast_d\tau=\tau$
then define $N_\tau:=\Delta_d^{-1}(C\tau+D\Delta_d)$ and you find that 
$M_j \Omega_\tau\,=\,\Omega_\tau N_\tau$.

That this fixed point set is a copy of $\HH_1$ in $\HH_2$ follows easily from
\cite[Hilfsatz III, 5.12, p.196]{Freitag}, but for convenience of the reader we provide an explicit description.

First of all, as in \cite[p.11]{HKW}, we introduce the matrix $4\times 4$ matrix $R_d$, 
where $I$ is the $2\times 2$ identity matrix
$$
R_d\,:=\,\begin{pmatrix} I&0\\ 0&\Delta_{d}\end{pmatrix},\qquad
\mbox{then}\quad E_d\,=\,R_dE_1{}^t\!R_d~.
$$
As $M_jE_d{}^t\!M_j=E_d$, we get $R_d^{-1}M_jR_d\in Sp(4,\RR)$, 
the (standard) real symplectic group of the (standard) alternating form $E_1$ and it is the matrix with rows $A,B\Delta_d$, $\Delta_d^{-1}C,\Delta^{-1}_dD\Delta_d$. 
Then, for the standard action $\ast_1$ of $Sp(4,\RR)$ on $\HH_2$, we find: 
$$
(R_d^{-1}M_jR_d)\ast_1\tau\,=\,(A\tau+B\Delta_d)(\Delta_d^{-1}C\tau+\Delta_d^{-1}D\Delta_d)^{-1}\,=\,
(A\tau+B\Delta_d)(C\tau+D\Delta_d)^{-1}\Delta_d~,
$$
hence $(R_d^{-1}M_jR_d)\ast_1\tau=M_j\ast_d \tau$.
So we need to describe the fixed points of $R_d^{-1}M_jR_d$ for the standard action on $\HH_2$. 
Let
$$
S_j\,:=\,\begin{pmatrix}(\sqrt{d}+1)/2&  \sqrt{d} &  0   &1\\
(\sqrt{d}-d)/2&-d&0&-\sqrt{d}\\
(\sqrt{d}-1)/2&-(\sqrt{d}+1)/2&1&-1\\
1&(1/\sqrt{d}+1)/2&1/\sqrt{d}&0
\end{pmatrix},\quad
\begin{pmatrix}
1/2&0&0&1\\
\sqrt{d}/2&0&0&-\sqrt{d}\\
0&-1/2&1&0\\
0&1/(2\sqrt{d})&1/\sqrt{d}&0
\end{pmatrix}~,
$$
for $j=3,4$ respectively.
One verifies that $S_3,S_4\in Sp(4,\RR)$ and that 
$M_j':=S_j^{-1}R_d^{-1}M_jR_dS_j$ is:
$$
M_3'\,=\,
\begin{pmatrix}0&1&0&0\\ 
-1&-1&0&0\\
0&0&-1&1\\
0&0&-1&0\end{pmatrix},\quad
M_4'\,=\,\begin{pmatrix}
0&1&0&0\\ 
-1&0&0&0\\
0&0&0&1\\
0&0&-1&0\end{pmatrix}~.
$$
Thus $M'_j\in Sp(2g,\RR)$ and $(R_d^{-1}M_jR_d)\ast_1\tau=\tau$ if and only if $M_j'(S_j^{-1}\ast_1\tau)=S_j^{-1}\ast_1\tau$. Thus $S_j$ induces a biholomorphic map between the fixed point sets of $M_j'$ and $R_d^{-1}M_jR_d$.  
As the matrix of $M_j'$ has blocks $D={}^t\!A^{-1}$, $B=C=0$, the fixed point set 
$\HH_2^{M'_j}$ of $M_j'$ is defined by $A\tau D^{-1}=\tau$, i.e.\ by 
$A\tau=\tau D$. 
It is then easy to find these fixed point sets: 
$$
\HH_2^{M_3'}\,=\,
\left \{\;\begin{pmatrix}2&-1\\-1&2\end{pmatrix}\tau\,:\;\tau\in\HH_1\right\}~,\qquad
\HH_2^{M_4'}\,=\,
\left \{\;\begin{pmatrix}1&0\\0&1\end{pmatrix}\tau\,:\;\tau\in\HH_1\right\}~,
$$
and one finds that they are indeed biholomorphic to $\HH_1$. \qed

\subsection{Polarizations and automorphisms} 
Recall that for a complex torus $A=\CC^g/\Lambda$ we can identify $\CC^g=\Lambda_\RR:=\Lambda\otimes_\ZZ\RR$. 
The scalar multiplication by $i=\sqrt{-1}$ on $\CC^g$ 
induces an $\RR$-linear map $J$ on $\Lambda_\RR$ with $J^2=-1$. 
An endomorphism of $A$ corresponds to a $\CC$-linear map $M$ on $\CC^g$ 
such that $M\Lambda\subset\Lambda$, equivalently,
after choosing a $\ZZ$-basis for $\Lambda$:
$$
\mbox{End}(A)\,=\,\{M\in M_{2g}(\ZZ):\; JM=MJ\,\}~,
$$
where $M_{2g}(\ZZ)$ is the algebra of $2g\times 2g$ matrices with integer coefficients. 

The N\'eron Severi group of $A$, a subgroup of $H^2(A,\ZZ)=\wedge^2H^1(A,\ZZ)=\wedge^2Hom(\Lambda,\ZZ)$, can be described similarly:
$$
\mbox{NS}(A)\,:=\,\{F\in M_{2g}(\ZZ):\, {}^tF\,=\,-F,\quad JF\,{}^t\!J\,=\,F\,\}~,
$$
where the alternating matrix $F\in \mbox{NS}(A)$ defines the bilinear form $(x,y)\mapsto xF\,{}^t\!y$. 
Moreover, $F$ is a polarization, 
i.e.\ the first Chern class of an ample line bundle, 
if $F\,{}^t\!J$ is a positive definite matrix.
In particular, $F$ is then invertible (in $M_{2g}(\QQ)$). 

It is now elementary to verify that if $E,F\in \mbox{NS}(A)$ and 
$E$ is invertible in $M_{2g}(\QQ)$, 
then $E^{-1}F\in \mbox{End}(A)_\QQ$ 
(cf.\ \cite[Proposition 5.2.1a] {BL} for an intrinsic description).
This result will be used in the proof of Theorem \ref{enddef}.

\

In Theorem \ref{enddef} we show that if $\tau\in\HH_{j,d}$ then the Abelian surface
$\mbox{End}(A_{\tau,d})_\QQ$ 
contains a quaternion  algebra (and not just the field $\QQ(\zeta_j)$!).
This is of course well-known (see for example the Exercise 4 of Section 9.4 of \cite{BL}), but we can also determine this quaternion algebra explicitly.
It allows us to find infinitely many families of $(1,d)$-polarized Abelian surfaces whose generic member is simple and whose endomorphism ring is an (explicitly determined) order in a quaternion algebra. To find the endomorphisms, we study first the N\'eron Severi group. 
Notice that in the proof of Theorem \ref{enddef} we don't need to know 
the period matrices of the deformations explicitly.

\subsection{Theorem} \label{enddef}
Let $j\in\{ 3,4\}$ and let $\tau\in  \HH_{j,d}$, so that the Abelian surface $A_{\tau,d}$
has an automorphism $\phi_j$ induced by $M_j$ (see Proposition \ref{defm}).

Then the endomorphism algebra of $A_{\tau,d}$ 
also contains an element $\psi_j$ with $\psi_j^2=d$ and such that
$\phi_j\psi_j=-\psi_j\phi_j$.
Moreover, for a general  $\tau\in  \HH_{j,d}$ one has
$$
End(A_{\tau,d})\,=\,\ZZ[\phi_j,\,\psi_j],\qquad
End(A_{\tau,d})_\QQ\,\cong\,(-j,\,d)_\QQ~,
$$
where 
$(a,b)_\QQ:=\QQ {\bf 1}\oplus\QQ {\bf i}\oplus\QQ {\bf j}
\oplus\QQ {\bf ij}$ is the quaternion algebra with
${\bf i}^2=a$, ${\bf j}^2=b$ and ${\bf ij}=-{\bf ji}$.

\ts
The N\'eron Severi group of an Abelian surface $A$ can also be described as
$$
\mbox{NS}(A)\,\stackrel{\cong}{\longrightarrow}\,
H^2(A,\ZZ)\cap H^{1,1}(A)\,\stackrel{\cong}{\longrightarrow}\,
\{\omega\in H^2(A,\ZZ)\,:\;(\omega,\omega_A^{2,0})\,=\,0\,\}~,
$$
where $(-,-)$ denotes the $\CC$-linear extension to $H^2(A,\CC)$  
of the intersection form on $H^2(A,\ZZ)$ and
we fixed a holomorphic $2$-form on $A$ so that
$H^{2,0}(A)=\CC \omega_A^{2,0}$.

The intersection form is invariant under automorphisms of $A$, so $(\phi_j^* x,\phi_j^*y)=(x,y)$ for all $x,y\in H^2(A,\ZZ)$, where $A=A_{\tau,d}$.
Moreover, by construction of $\phi_j$, we have that 
$\phi_j^*\omega^{2,0}_A=\omega^{2,0}_A$, so 
$\omega^{2,0}_A\in H^2(A,\CC)^{\phi_j^*}$, 
the subspace of $\phi_j$-invariant classes.
Therefore any integral class which is perpendicular to the $\phi_j$-invariant classes 
is in particular perpendicular to $\omega_A^{2,0}$ 
and thus must be in $\mbox{NS}(A)$:
$$
\left(H^2(A,\ZZ)^{\phi_j^*}\right)^\perp\,:=\,
\{\omega\in H^2(A,\ZZ):\;(\omega,\theta)\,=\,0,\quad
\forall\,\theta\,\in H^2(A,\ZZ)\;\mbox{with}\;\phi_j^*\theta\,=\,\theta\}\,\subset\,
\mbox{NS}(A)~.
$$
The eigenvalues of $\phi_j^*$ on $H^1(A,\CC)=H^{1,0}(A)\oplus \overline{H^{1,0}(A)}$ are $\zeta_j$ and $\zeta_j^{-1}$, both with multiplicity two. Thus the eigenvalues of $\phi^*$ on $H^2(A,\CC)=\wedge^2H^1(A,\CC)$ are
$\zeta_j^2,\zeta_j^{-2}$, with multiplicity one, and $1$ with multiplicity $4$. In particular
$\left(H^2(A,\ZZ)^{\phi_j^*}\right)^\perp$ is a free $\ZZ$-module of rank $2$, it is 
the kernel in $H^2(A,\ZZ)$ of $(\phi_3^*)^2+\phi_3^*+1$ in case $j=3$ and of $(\phi^*_4)^2+1$ in case $j=4$. 
Identifying $H^2(A,\ZZ)$ with the alternating bilinear $\ZZ$-valued maps on $\Lambda_j\cong \ZZ^4$, the action of $\phi^*$ is given by $M_j\cdot F:=M_jF{}^tM_j$,
where $F$ is an alternating $4\times 4$ matrix with integral coefficients. 
It is now easy to find a basis $E_{j,1}$, $E_{j,2}$ of the $\ZZ$-module 
$\left(H^2(A,\ZZ)^{\phi_j^*}\right)^\perp$. 
Since $E_{d}$ defines a polarization on $A$,
the matrices $E_{d}^{-1}E_{j,k}$, $k=1,2$, 
are the images under $\rho_r$ of elements in $End(A)_\QQ$ 
(cf.\ \cite[Proposition 5.2.1a] {BL}).
In this way we found that for any $\tau\in \HH_{j,d}$, the Abelian surface $A=A_{\tau,d}$
has an  endomorphism $\psi_j$ defined by 
the matrix $\rho_r(\psi_j)$ below:
$$
\rho_r(\psi_3)\,=\,
\begin{pmatrix}0&d&0&0\\1&0&0&0\\
0&0&0&d\\0&0&1&0\end{pmatrix},\qquad
\rho_r(\psi_4)\,=\,\begin{pmatrix}0&0&0&-d\\ 0&0&1&0\\
0&d&0&0\\-1&0&0&0\end{pmatrix}
~.
$$
It is easy to check that $\rho_r(\psi_j)^2=d$ and that $M_4\rho_r(\psi_4)=-\rho_r(\psi_4)M_4$, whereas $(1+2M_3)\rho_r(\psi_3)=-\rho_r(\psi_3)(1+2M_3)$ (and notice that $(1+2M_3)^2=-3$). Therefore $(-j,d)_\QQ\subset End(A)_\QQ$
(in fact, $M_4^2=-1$, but $(-1,d)_\QQ\cong (-4,d)_\QQ$).
As $(-j,d)_\QQ$ is a (totally) indefinite quaternion algebra (so of type $II$), for  general
$\tau\in\HH_{j,d}$ the Abelian surface $A=A_{\tau,d}$  has 
$(-j,d)_\QQ= \mbox{End}(A)_\QQ$ by \cite[Theorem 9.9.1]{BL}.
Therefore if $\phi\in \mbox{End}(A)$, then $\rho_r(\phi)$ is both a matrix with integer coefficients and it is a linear combination of 
$I$, $M_j=\rho_r(\phi_j)$, $\rho_r(\psi_j)$ and $M_j\rho_r(\psi_j)$ with rational coefficients. It is then
easy to check that $\mbox{End}(A)$ is as stated in Theorem \ref{enddef}.
\qed

\subsection{A table}\label{cordisc}
Using Magma (\cite{magma}), we found that for  the following $d\leq 20$
the quaternion algebras $(-1,d)_\QQ$ and $(-3,d)_\QQ$ are  skew fields:
$$
\left|
\begin{array}{cc}
d&\mbox{discriminant}(-1,d)_\QQ\\
3,6,15&  6\\
7,14&14\\
11&22\\
19&38\\
\end{array}\right|~, \qquad\qquad
\left|
\begin{array}{cc}
d&\mbox{discriminant}(-3,d)_\QQ\\
2,6,8,14,18&  6\\
5,15,20&15\\
10&10\\
11&33\\
17&51\\
\end{array}\right|~.
$$
Moreover, for $d\leq 20$, $\mbox{End}(A)$ is never a maximal order in $(-1,d)_\QQ$, and it is a maximal order in $(-3,d)_\QQ$ if and only if $d=2,5,11,17$.

In particular, for $\tau\in \HH_{3,2}$ the Abelian surface $A_{\tau,2}$ has a $(1,2)$-polarization invariant by an automorphism of order three induced by $M_3$ and  $\mbox{End}(A_{\tau,2})=\cO_6$, the maximal order in the quaternion algebra with discriminant 6, for general $\tau\in\HH_{3,2}$.  
After a discussion of an equivariant map $\overline{\psi}_D$ 
of a moduli space of Abelian surfaces to a projective space, 
we will describe the image of $\HH_{3,2}$ in Section \ref{projmodel}.

\

\section{The level moduli space}

\subsection{The moduli space of (1, d)-polarized Abelian surfaces}\label{mod1d}
The integral symplectic group with respect to $E_{d}$ is defined as
$$
\tilde{\Gamma}^0_d
\,:=\,\{\,M\in\,GL(4,\ZZ)\,:\,ME_{d}{}^t\!M\,=\,E_{d}\,\}~.
$$
This group acts on the Siegel space by (\cite[Equation (1.4)]{HKW}):
$$
\tilde{\Gamma}^0_d\times\HH_2\,\longrightarrow\,\HH_2,\qquad
\begin{pmatrix}A&B\\C&D\end{pmatrix}\ast_d\tau\,:=\,
(A\tau+B\Delta_d)(C\tau+D\Delta_d)^{-1}\Delta_d~.
$$
Notice that for $d=1$ one finds the standard action of the symplectic group on $\HH_2$.
The quotient space (in general a singular quasi-projective 3-dimensional algebraic variety) is the moduli space $\cA_d^0$ of pairs $(A,H)$ where $A$ is an Abelian surface and $H$ is a polarization of type $(1,d)$ (\cite[Theorem 1.10(i)]{HKW}).

For the study of this moduli space, and of certain `level' covers of it,
we use the standard action of $Sp(4,\RR)$ on $\HH_2$
which is $\ast_1$.
For this, as in the proof of Proposition \ref{defm} (cf.\  \cite[p.11]{HKW}),
we use the matrix $4\times 4$ matrix $R_d$. 
Then $\Gamma_{1,d}^0:=R_d^{-1}\tilde{\Gamma}^0_dR_d\in Sp(4,\RR)$
is a subgroup of the (standard) real symplectic group of the 
(standard) alternating form $E_1$ and
we have $(R_d^{-1}MR_d)\ast_1\tau=M\ast_d \tau$ for all $M\in \tilde{\Gamma}^0_d$.
Therefore
$$
\cA_d^0\,:=\,\tilde{\Gamma}^0_d\backslash\HH_2\,\cong\,
\Gamma^0_{1,d}\backslash\HH_2~,
$$
where the actions are $\ast_d$ and $\ast_1$ respectively.

\subsection{Congruence subgroups} 
We now follow \cite{BL} for the definition of coverings of the moduli space and 
maps to projective space. 
Recall that we defined a group $\tilde{\Gamma}^0_d$ in Section \ref{mod1d}
of matrices with integral coefficients which preserve the alternating form $E_d$.
We will actually be interested in the form $2E_2$, which is preserved by the same group.
With the notation from \cite[8.1, p.212]{BL} we thus have:
$$
\tilde{\Gamma}^0_2\,=\,\Gamma_D\,=\,Sp_4^D(\ZZ),\qquad D\,=\,diag(2,4)\,=\,2\Delta_2~.
$$
It is easy to check that
$$
\ZZ^4\tilde{D}^{-1}\,=\,\{x\in\QQ^4: x(2E_2)y\in\ZZ,\;\forall\;y\,\in\,\ZZ^4\,\},
\qquad\tilde{D}\,:=\,\begin{pmatrix}D&0\\0&D\end{pmatrix}~.
$$
Let $T(2,4)$ be the following quotient of $\ZZ^4$: 
$$
T(2,4)\,=\,(\ZZ^4\tilde{D}^{-1})/\ZZ^4\,\cong(\ZZ/2\ZZ\times\ZZ/4\ZZ)^2,
$$
The group $\Gamma_D$ acts on this quotient and we define
$$
\Gamma_D(D)\,:=\,\ker(\Gamma_D\,\longrightarrow\,Aut(T(2,4)))~.
$$
One verifies easily that
{\renewcommand{\arraystretch}{1.7}
$$
\begin{array}{rcl}
\Gamma_D(D)&=&
\{M\in\Gamma_D:\, \tilde{D}^{-1}M\equiv \tilde{D}^{-1}\;\mbox{mod}\;M_4(\ZZ)\,\}
\\
&=& 
\left\{\,M\,=\,\begin{pmatrix}I+D\alpha&D\beta\\D\gamma&I+D\delta\end{pmatrix}
\,\in\,\Gamma_D:\;\alpha,\beta,\gamma,\delta\in M_2(\ZZ)\,\right\}~.
\end{array}
$$
} 
This shows that $\Gamma_D(D)$ is the subgroup as defined in \cite[Section 8.3]{BL} 
(see also \cite[Section  8.8]{BL}). 
The alternating form $E_2$ defines a `symplectic' form $<,>$ on $T(2,4)$ with values in the fourth-roots of unity (cf.\ \cite[Section 3.1]{Barth}). 
For this we write (cf.\ \cite[Section 2.1]{Barth})
$$
T(2,4)\,=\,K\times \hat{K},\qquad
K\,=\,\ZZ/2\ZZ\,\times\, \ZZ/4\ZZ,\quad
\hat{K}\,=\,\mbox{Hom}(K,\CC^\ast)\,\cong\,\ZZ/2\ZZ\,\times\, \ZZ/4\ZZ~,
$$
and the symplectic form is
$$
<-,->:\,T(2,4)\,\times\,T(2,4)\,\longrightarrow\,\CC^\ast,\qquad
<(\sigma,l),(\sigma',l')>\,:=\,l'(\sigma)l(\sigma')^{-1}~.
$$
We denote by $Sp(T(2,4))$ the subgroup of $Aut(T(2,4))$ of automorphisms which preserve
this form.

\subsection{Lemma}\label{redu24}
The reduction homomorphism 
$$
\Gamma_D\,\longrightarrow\,Sp(T(2,4))
$$
is surjective. 
Hence $\Gamma_D/\Gamma_D(D)\cong Sp(T(2,4))$,
this is a finite group of order $2^93^2$.

\ts
As the symplectic form is induced by $E_2$, we have $\im(\Gamma_D)\subset Sp(T(2,4))$.
In \cite[Proposition 3.1]{Barth} generators $\phi_i$, $i=1,\ldots,5$ of $Sp(T(2,4))$ are given. It is easy to check that the following matrices are in $G_D$
and induce these automorphisms on $T(2,4)$:
$$
\left(\begin{smallmatrix}
1&0&0&0\\0&0&0&-1\\0&0&1&0\\0&1&0&0\end{smallmatrix}
\right),\quad
\left(\begin{smallmatrix}1&0&0&0\\0&1&0&0\\0&0&1&0\\0&1&0&1\end{smallmatrix}
\right),\quad
\left(\begin{smallmatrix}0&0&-1&0\\0&1&0&0\\1&0&0&0\\0&0&0&1\end{smallmatrix}
\right),\quad
\left(\begin{smallmatrix}1&0&0&0\\0&1&0&0\\1&0&1&0\\0&0&0&1\end{smallmatrix}
\right),\quad
\left(\begin{smallmatrix}1&0&0&0\\-2&1&0&0\\0&0&1&1\\0&0&0&1\end{smallmatrix}
\right)~.
$$
The order of $Sp(T(2,4))$ is determined in \cite[Proposition 3.1]{Barth}.
\qed

\subsection{The subgroup $\Gamma_D(D)_0$}
We define a normal subgroup of $\Gamma_D(D)$ by:
$$
\Gamma_D(D)_0\,:=\,\ker(\phi:\Gamma_D(D)\,\longrightarrow\,(\ZZ/2\ZZ)^4),\qquad
\phi(M)\,=\,(\beta_0,\gamma_0)\,:=\,(\beta_{11},\beta_{22},\gamma_{11},\gamma_{22})~,
$$
where $M\in \Gamma_D(D)$ is as above.
Since $D$ has even coefficients, $D=2diag(1,2)$, it is easy to check that $ \phi$ is a homomorphism. Moreover, $\phi$ is surjective since the matrix with $\alpha=\gamma=\delta=0$ and $\beta=diag(a,b)$ ($a,b\in\ZZ$) 
is in $\Gamma_D(D)$ and maps to $(a,b,0,0)$, similarly the matrix with $\alpha=\beta=\delta=0$ and $\gamma=diag(a,b)$ is also in $\Gamma_D(D)$ 
and maps to $(0,0,a,b)$. 
It follows that $\Gamma_D(D)/\Gamma_D(D)_0\cong(\ZZ/2\ZZ)^4$.

The groups $\Gamma_D,\Gamma_D(D)$ and $\Gamma_D(D)_0$ are denoted by 
$G_\ZZ,G_\ZZ(e)$ and $G_\ZZ(e,2e)$ in \cite[V.2, p.177]{IgusaT}.
In \cite[V.2 Lemma 4]{IgusaT} one finds that $\Gamma_D(D)_0$
is in fact a normal subgroup of $\Gamma_D$.
There is an exact sequence of groups:
$$
0\,\longrightarrow\,\Gamma_D(D)/\Gamma_D(D)_0\,\longrightarrow\,
\Gamma_D/\Gamma_D(D)_0\,\longrightarrow\,
\Gamma_D/\Gamma_D(D)\,\longrightarrow\,0~.
$$

The group $\Gamma_D$ act on $\HH_2$ in a natural way but to get the standard action $\ast_1$ one must conjugate these groups by a matrix $R_D$ with diagonal blocks $I$, $D$ and one obtains
the groups
$$
G_D\,=\,R_D^{-1}\Gamma_DR_D,\qquad 
G_D(D)\,=\,R_D^{-1}\Gamma_D(D)R_D,\qquad 
G_D(D)_0\,=\,R_D^{-1}\Gamma_D(D)_0R_D~,
$$
see \cite[Section 8.8, 8.9]{BL}. 

The main result from \cite[section 8.9]{BL} is Lemma 8.9.2 which asserts that the holomorphic map given by theta-null values
$$
\psi_D:\HH_2\,\longrightarrow\,\PP^7,\qquad
\tau\,\longrightarrow\,\left(\ldots:\,\vartheta[{}^l_0](0,\tau)\,:\ldots\right)_{l\in K}~,
$$
where $l$ runs over $K=D^{-1}\ZZ^2/\ZZ^2$
and where the theta functions $\vartheta[{}^l_0](v,\tau)$ are defined in 
\cite[8.5, Formula (1)]{BL}, factors over a holomorphic map
$$
\overline{\psi}_D:\, \cA_D(D)_0\,:=\,\HH_2/\Gamma_D(D)_0\,\cong\,
\HH_2/G_D(D)_0\,\longrightarrow\,\PP^7~.
$$

\subsection{Group actions}
The finite group $\Gamma_D/\Gamma_D(D)_0$ acts on $\cA_D(D)_0$.
On $\PP^7$ the Heisenberg group $\ccH(D)$, a non-Abelian central 
extension of $T(2,4)$ by $\CC^*$, acts (\cite[Section 6.6]{BL}). 
This action is induced by irreducible representation
(called the Schr\"odinger representation) of $\ccH(D)$
on the vector space $V(2,4)$ of complex valued functions on the subgroup $K$ 
of $T(2,4)$ (\cite[Section 6.7]{BL})
$$
\rho_D:\,\ccH(D)\,\longrightarrow\,GL(V(2,4))~.
$$
In \cite[Section 2.1]{Barth}) the action of generators of $\ccH(D)$ on $\PP V(2,4)=\PP^7$ are given explicitly, also the linear map $\tilde{\iota}\in GL(V(2,4))$ which sends the delta functions $\delta_l\mapsto \delta_{-l}$ ($l\in K$) is introduced there (cf.\ Sections \ref{m24}, \ref{hact}).

The normalizer of the Heisenberg group (in the Schr\"odinger representation)
is by definition the group
$$
N(\ccH(D))\,:=\,\{\,\gamma\in Aut(\PP V(2,4)):\;\gamma \rho_D(\ccH(D))\gamma^{-1}
\,\subset\,\rho_D(\ccH(D))\,\}~.
$$
The group $N(\ccH(D))$ maps onto $Sp(T(2,4))$ with kernel isomorphic to $T(2,4)$.
The elements in this kernel are obtained as interior automorphisms:
$\gamma=\rho_D(h)$, for some $h\in \ccH(D)$. 
Explicit generators of $N(\ccH(D))$ are given in \cite[Table 8]{Barth} 
(but there seem to be some misprints in the action of the generators on $\ccH(D)$ 
in the lower left corner of that table).
Let $N(\ccH(D))_2$ be the subgroup of 
$N(\ccH(D))$ of elements which commute with $\tilde{\iota}$.
The group $N(\ccH(D))_2$ is an extension of $Sp(T(2,4))$ by the 2-torsion subgroup (isomorphic to $(\ZZ/2\ZZ)^4$) of
$T(2,4)$ and  $\sharp N(\ccH(D))_2=2^{13}3^2$.

We need the following result.

\subsection{Proposition}\label{equi} There is an isomorphism 
$\gamma:G_D/G_D(D)_0\,\cong\,N(\ccH(D))_2$, $M'\mapsto \gamma_{M'}$ such that the 
map $\overline{\psi}_D$ is equivariant for the action of these groups.
So if we denote by $\tilde{\gamma}$ the composition
$$
\tilde{\gamma}:\,\Gamma_D/\Gamma_D(D)_0\,
\stackrel{\cong}{\longrightarrow}\,G_D/G_D(D)_0\,
\stackrel{\gamma}{\longrightarrow}\,N(\ccH(D))~,
$$
then $\overline{\psi}_D(M\ast \tau)=\tilde{\gamma}_M\overline{\psi_D}(\tau)$ where
$\ast$ denotes the action of $\Gamma(D)$ on $\HH_2$.

\ts
Let $\cL_\tau=L(H,\chi_0)$ be the line bundle on 
$A_{\tau,2}:=\CC^2/(\ZZ^4\Omega_\tau)$
which has Hermitian form $H$ with $E_2=\mbox{Im} H$ 
(so it defines a polarization of type $(1,2)$) and the quasi-character $\chi_0$ is as in  
\cite[3.1, Formula (3)]{BL} for the decomposition 
$\Lambda=\ZZ^2\tau\oplus \ZZ^2\Delta_2$.
According to \cite[Remark 8.5.3d]{BL}, the theta functions 
$\vartheta[{}^l_0](v,\tau)$
are a basis of the vector space of classical theta functions for the line bundle 
$\cL^{\otimes 2}_\tau$. As $\chi_0$ takes values in 
$\{\pm 1\}$ one has $\cL_\tau^{\otimes 2}=L(2H,\chi_0^2=1)$, so it
is the unique line bundle with first Chern class $2E_2$ 
and trivial quasi-character. 
Thus if $M\in G_D$ and $\tau'=M\ast_1\tau$ then 
$\phi_M^*\cL^{\otimes 2}_\tau\cong\cL^{\otimes 2}_{\tau'}$, 
where $\phi_M:A_{\tau',2}\rightarrow A_{\tau,2}$ 
is the isomorphism defined by $M$.
Notice that $\cL_\tau$ and $\cL_{\tau}^{\otimes 2}$ are symmetric line bundles
(\cite[Corollary 2.3.7]{BL}).

Let $\cG(\cL_\tau^{\otimes 2})$ be the theta group (\cite[Section 6.1]{BL}), it has an
irreducible linear representation $\tilde{\rho}$  
on $H^0(A_{\tau,2},\cL_\tau^{\otimes 2})$ (\cite[Section 6.4]{BL}).

A theta structure $b:\cG(\cL_\tau^{\otimes 2})\rightarrow \ccH(D)$ is an isomorphism
of groups which is the identity on their subgroups $\CC^\ast$.
A theta structure $b$ defines an
isomorphism $\beta_b$, unique up to scalar multiple (\cite[Section 6.7]{BL}), which intertwines
the actions of $\cG(\cL^{\otimes 2})$ and $\ccH(D)$:
$$
\beta_b:H^0(A_{\tau,2},\cL_\tau^{\otimes 2})\,\longrightarrow \,V(2,4),\qquad
\beta_b\tilde{\rho}(g)\,=\,\rho_D(b(g))\beta_b
\quad(\forall g\in \cG(\cL^{\otimes 2}_\tau))~.
$$ 
A symmetric theta structure  (\cite[Section 6.9]{BL}) is a theta structure
which is compatible with the action of $(-1)\in \mbox{End}(A_{\tau,2})$ 
on the symmetric line bundle $\cL_\tau^{\otimes 2}$ and the map 
$\tilde{\iota}\in GL(V(2,4))$ defined in \cite[Section 2.1]{Barth}.

For $\tau\in\HH_2$, define an isomorphism 
$\beta_\tau:H^0(A_{\tau,2},\cL_\tau^{\otimes 2})\,\rightarrow \,V(2,4)$
by sending the basis vectors $\vartheta[{}^l_0](v,\tau)$ to the delta functions $\delta_l$ for $l\in K$.
From the explicit transformation formulas for the theta functions under translations by points in $A_{\tau,2}$ one finds that for $g\in \cG(\cL^{\otimes 2}_\tau)$ the map
$\beta_\tau\tilde{\rho}(g)\beta_\tau^{-1}$ 
acts as an element, which we denote by $b_\tau(g)$, 
of the Heisenberg group $\ccH(D)$ acting on 
$V(2,4)$.
This map $b=b_\tau: \cG(\cL^{\otimes 2}_\tau)\rightarrow\ccH(D)$ is a theta structure and $\beta_\tau\tilde{\rho}(g)\,=\,\rho_D(b_\tau(g))\beta_\tau$, moreover it is symmetric since 
$\theta[{}^l_0](-v,\tau)=\theta[{}^{-l}_0](v,\tau)$.

For $M\in G_D$ and $\tau'=M\ast_1\tau$ we have an isomorphism $\beta_{\tau'}$
and the composition $\gamma_M:=\beta_{\tau'}\phi_M^*\beta_{\tau}^{-1}\in GL(V(2,4))$,
is an element of $N(\ccH)$ since $\phi_M^*$ induces an isomorphism 
$\cG(\cL^{\otimes 2}_\tau)\rightarrow\cG(\cL^{\otimes 2}_{\tau'})$.
In fact $\gamma_M\in N(\ccH)_2$ since the theta structures $\beta_{\tau}$, $\beta_{\tau'}$
are symmetric and $\phi_M$ commutes with $(-1)$ on the abelian varieties.

From \cite[Proposition 6.9.4]{BL} it follows that the group generated by the $\gamma_M$
is contained in an extension of $Sp(T(2,4))$ by $(\ZZ/2\ZZ)^4$. 
The map $M\mapsto \gamma_M\in Aut(\PP(V(2,4))$ is thus a (projective) representation of
$G_D$  whose image is contained in $N(\ccH)_2$ and which, by construction, is equivariant
for $\overline{\psi}_D$.
Unwinding the various definitions, we have shown that $\gamma_M$ maps the point
$(\ldots:\theta[{}^l_0](v,\tau):\ldots)$ to the point 
$(\ldots:\theta[{}^l_0]({}^t\!(C\tau+D)v,M\ast_1\tau):\ldots)$ where $M$ has block form $A,\ldots,D$.
From the classical theory of transformations of theta functions (as in \cite[Section 8.6]{BL})
one now deduces that $M\mapsto \gamma_M$ provides the desired isomorphism of groups.
Notice that the element $-I\in G_D$, which acts trivially on $\HH_2$,
maps to $\tilde{\iota}\in N(\ccH(D))_2$ which acts trivially on the subspace $\PP^5\subset\PP^7$ of even theta functions.
\qed

\

\section{A projective model of a Shimura curve}\label{projmodel}

\subsection{Barth's variety $M_{2,4}$}\label{m24}
We choose projective coordinates $x_1,\ldots,x_8$ on 
$\PP^7=\PP V(2,4)$ as in \cite[\S 2.1]{Barth}. 
The map $\tilde{\iota}\in Aut(\PP^7)$ is then given by
$$
\tilde{\iota}(x)\,=\,(\,x_1:\,x_2:\,x_3:\,x_4:\,x_5:\,x_6:\,-x_7:\,-x_8)~.
$$
It has two eigenspaces which correspond to the even and odd theta functions. The image of $\overline{\psi}_D$ lies in the subspace
$\PP^5=\PP V(2,4)_+$ of even functions which is defined by $x_7=x_8=0$.
We use $x_1,\ldots,x_6$ as coordinates on this $\PP^5$.
Let
$$
f_1\,:=\,- x_1^2x_2^2\,+\,x_3^2x_4^2\,+\,x_5^2x_6^2,\qquad
f_2\,:=\,-(x_1^4+x_2^4)\,+\,x_3^4+x_4^4+x_5^4+x_6^4~.
$$
Then Barth's variety of theta-null values is defined as (\cite[(3.9)]{Barth})
$$
M_{2,4}\,:=\,\{x\in\PP^5:\quad f_1(x)\,=\,f_2(x)\,=\,0\,\}~.
$$
The image of $\overline{\psi}_D(\HH_2)$ is a quasi-projective variety and the closure of its image is $M_{2,4}$.

\subsection{The Heisenberg group action}\label{hact}
Recall that $T(2,4)=\ZZ^4\tilde{D}^{-1}/\ZZ^4$ and let 
$\sigma_1,\sigma_2,\tau_1,\tau_2\in T(2,4)$ be the 
images of $e_1/2,e_2/4,e_3/2,e_4/4$.
We denote certain lifts of the generators
$\sigma_1,\ldots,\tau_2$ of $T(2,4)$ to $\ccH(D)$
by $\tilde{\sigma}_1,\ldots,\tilde{\tau}_2$.
These lifts act, in the Schr\"odinger representation, on $\PP^7=\PP V(2,4)$ as follows
(see \cite[Table 1]{Barth}):
{\renewcommand{\arraystretch}{1.3}
$$
\begin{array}{rclrrrrrrrr}
\tilde{\sigma}_1(x)&=&(&x_2:&x_1:&x_4:&x_3:&x_6:&x_5:&x_8:&x_7)~,\\
\tilde{\sigma}_2(x)&=&(&x_3:&x_4:&x_1&x_2:&x_7:&x_8:&-x_5:&-x_6)~,\\
\tilde{\tau}_1(x)&=&(&x_1:&-x_2:&x_3:&-x_4:&x_5:&-x_6:&x_7:&-x_8)~,\\
\tilde{\tau}_2(x)&=&(&x_5:&x_6:&ix_7:&ix_8:&x_1:&x_2:&ix_3:&ix_4)~,\\
\end{array}
$$
}
where $x=(x_1:\ldots:x_8)\in\PP^7$ and $i^2=-1$.
For any $g=(a,b,c,d)\in T(2,4)$ one then finds the action of a lift $\tilde{g}$ of $g$ 
by defining $\tilde{g}:=\tilde{\sigma}_1^a\cdots\tilde{\tau}_2^d$.

\subsection{Proposition}\label{mu3}
Let $\tilde{\mu}_3$ on $\PP^7$ be the projective transformation
defined as $\tilde{\mu}_3:x\mapsto$
$$
\big(x_3-ix_4:\,x_3+ix_4:\,
\zeta x_5-\zeta^3 x_6:\,
\zeta x_5+\zeta^3 x_6:\,
x_1-ix_2:\,x_1+ix_2:\, 
\zeta^3 x_7+\zeta x_8:\,\zeta^3 x_7-\zeta x_8
\big)~,
$$
where $\zeta$ is a primitive $8$-th root of unity (so $\zeta^4=-1$) and 
$i:=\zeta^2$.
Then $\tilde{\mu}_3\in N(\ccH(D))_2$ and with $M_3$ as in Section \ref{prelc} we have
$$
\tilde{\gamma}_{M_3}\,=\,\tilde{h}\tilde{\mu}_3\tilde{h}^{-1}
$$ 
for some $\tilde{h}\in \ker(N(\ccH(D))_2\rightarrow Sp(T(2,4))$.

\ts
The map $M_3:\ZZ^4\rightarrow\ZZ^4$ from section \ref{prelc} induces the (symplectic) automorphism $\overline{M}_3$ of $T(2,4)$ given by
(recall that we used row vectors, so for example $e_4M_3=-e_2-e_4$ and thus
$\tau_2\mapsto -\sigma_2-\tau_2$):
$$
\sigma_1\,\longmapsto\,-\sigma_1-\tau_1,\quad
\sigma_2\,\longmapsto\,\tau_2,\quad
\tau_1\,\longmapsto\,\sigma_1,\quad
\tau_2\,\longmapsto\,-\sigma_2-\tau_2~.
$$
Now one verifies that, as maps on $\CC^8$, one has 
$$
\tilde{\mu}_3 \tilde{\sigma}_1\tilde{\mu}_3^{-1}=
i\tilde{\sigma}_1^{-1}\tilde{\tau}_1^{-1},
\quad
\tilde{\mu}_3 \tilde{\sigma}_2\tilde{\mu}_3^{-1}=\tilde{\tau}_2,\quad
\tilde{\mu}_3 \tilde{\tau}_1\tilde{\mu}_3^{-1}=\tilde{\sigma}_1,\quad
\tilde{\mu}_3 \tilde{\tau}_2\tilde{\mu}_3^{-1}
=\zeta\tilde{\sigma}_2^{-1}\tilde{\tau}_2^{-1}~.
$$

Hence $\tilde{\mu}_3\in Aut(\PP^7)$ is in the normalizer $N(\ccH)$
and it is a lift of $\overline{M}_3\in Sp(T(2,4))$.
One easily verifies that it commutes with the action of $\tilde{\iota}$ on $\PP^7$ so
$\tilde{\mu}_3\in N(\ccH)_2$.
Any other lift of $\overline{M}_3$ to $Aut(\PP^7)$ which commutes with $\tilde{\iota}$
is of the form $\tilde{g}\tilde{\mu}_3$ for some $g\in T(2,4)$ with $2g=0$.
Since $\overline{M}_3^2+\overline{M}_3+I=0$, the map $h\mapsto (\overline{M}_3+I)h$
is an isomorphism on the two-torsion points in $T(2,4)$. Thus there is an $h\in T(2,4)$, with $2h=0$, such that $g=(\overline{M}_3+I)h$. As $\tilde{\mu}_3\tilde{h}\tilde{\mu}_3^{-1}=\tilde{k}$, where $k=\overline{M}_3h$
and thus $k=g+h$, it follows that
$\tilde{h}\tilde{\mu}_3\tilde{h}^{-1}=\tilde{g}\tilde{\mu}_3$.
\qed

\subsection{Fixed points and eigenspaces}
The map $\overline{\psi}_D$ is equivariant for the actions of $\Gamma_D$ 
and $N(\ccH)_2$.
Hence the fixed points of $M_3$ in $\HH_2$, which parametrize abelian surfaces with quaternionic multiplication, map to the fixed points of $\tilde{\gamma}_{M_3}=\tilde{h}\tilde{\mu}_3\tilde{h}^{-1}$ in $\PP^7$.
Conjugating $M_3$ by an element $N\in \Gamma_D$ such that $\tilde{\gamma}_N=\tilde{h}$ (as in Proposition \ref{mu3}), 
we obtain an element of order three $M_3'\in \Gamma_D$ whose fixed point locus 
$\HH_2^{M_3'}$ also consists of
period matrices of Abelian surfaces with QM by $\cO_6$ and the image 
$\overline{\psi}_D(\HH_2^{M_3'})$ consists of fixed points of $\tilde{\mu}_3$.
The following lemma identifies this fixed point set.

\subsection{Theorem}\label{p1qm}
Let $\PP^1_{QM}\subset\PP^5$ be the projective line parametrized by
$$
\PP^1\stackrel{\cong}{\longrightarrow}\PP^1_{QM},\qquad
(x:y)\,\longmapsto\,p_{(x:y)}:=\,(\sqrt{2}x:\sqrt{2}y:x+y:i(x-y):x-iy:x+iy)~.
$$
Then $\PP^1_{QM}\subset M_{2,4}$ is a Shimura curve that 
parametrizes Abelian surfaces with QM by $\cO_6$, the maximal order in the quaternion algebra of discriminant $6$.

The following two elements $\tilde{\nu}_1,\tilde{\nu}_2\in N(\ccH(D))_2$,
$$
\tilde{\nu}_1(x)\,=\,(x_5+x_6,-x_5+x_6,\zeta (x_3-x_4),\zeta(x_3+x_4),
x_1+x_2,x_1-x_2, (\zeta(-x_7+x_8),\zeta(x_7+x_8))
$$
$$
\tilde{\nu}_2(x)\,=\,(x_4,-x_3,\zeta^3x_6,\zeta^3x_5,
i x_1,-i x_2,\zeta^3 x_7,\zeta^3 x_8)~,
$$
restrict to maps in $Aut(\PP^1_{QM})$ which generate a subgroup
isomorphic to the symmetric group $S_4\subset Aut(\PP^1_{QM})$.

\ts
The subspace $\PP^5$  is mapped into itself by $\tilde{\mu}_3$. 
The restriction $\mu_3$ of $\tilde{\mu}_3$
to $\PP^5$ has three eigenspaces on $\CC^6$, each 2-dimensional.
The eigenspace of $\mu_3$ with eigenvalue $\sqrt{2}:=\zeta+\zeta^7$ 
is the only eigenspace whose projectivization $\PP^1_{QM}$ is contained in $M_{2,4}$.
Thus $\overline{\psi}_D(\HH_2^{M_3'})\subset \PP^1_{QM}$ and we have equality 
since the locus of Abelian surfaces with QM by $\cO_6$ in $\cA_D(D)_0$ (in fact in any level moduli space) is known to be a compact Riemann surface.

The maps $\tilde{\nu}_1$, $\tilde{\nu}_2$
commute with $\tilde{\iota}$ and moreover: 
$$
\tilde{\nu}_1 \tilde{\sigma}_1\tilde{\nu}_1^{-1}=
-\tilde{\sigma}_1\tilde{\tau}_2^{2},
\quad
\tilde{\nu}_1 \tilde{\sigma}_2\tilde{\nu}_1^{-1}=
i\tilde{\sigma_1}\tilde{\sigma_2}^2\tilde{\tau}_2,
\quad
\tilde{\nu}_1 \tilde{\tau}_1\tilde{\nu}_1^{-1}=-\tilde{\sigma}_2^2\tilde{\tau_1},\quad
\tilde{\nu}_1 \tilde{\tau}_2\tilde{\nu}_1^{-1}
=\zeta\tilde{\sigma}_1\tilde{\sigma}_2^3\tilde{\tau}_1\tilde{\tau}_2~,
$$
$$
\tilde{\nu}_2 \tilde{\sigma}_1\tilde{\nu}_2^{-1}=
-\tilde{\tau}_1\tilde{\tau}_2^{2},
\quad
\tilde{\nu}_2 \tilde{\sigma}_2\tilde{\nu}_2^{-1}=
\zeta\tilde{\sigma_1}\tilde{\sigma_2}\tilde{\tau}_2,
\quad
\tilde{\nu}_2 \tilde{\tau}_1\tilde{\nu}_2^{-1}=
-\tilde{\sigma}_1\tilde{\sigma}_2^2\tilde{\tau}_2^2,\quad
\tilde{\nu}_2 \tilde{\tau}_2\tilde{\nu}_2^{-1}
=\tilde{\tau}_1\tilde{\tau}_2^3~,
$$
hence they are in $N(\ccH)_2$.
The maps $\nu_1$, $\nu_2$ have order $4$ and $3$ respectively in $Aut(\PP^7)$
and map $\PP^1_{QM}$ into itself, in fact,
the induced action on $\PP^1_{QM}$ is:
$$
\tilde{\nu}_i p_{(x:y)}\,=\,p_{\nu_i(x:y)}\quad
\mbox{with}\quad
\nu_1(x:y)\,:=\,(x:iy),\quad \nu_2(x:y)\,:=\,(i(x-y):-(x+y))~.
$$
We verified that $\nu_1,\nu_2\in Aut(\PP^1)$ generate a subgroup which is 
isomorphic to the symmetric group $S_4$ (to obtain this isomorphism, one may use the action of the $\nu_i$ on the four irreducible factors in $\QQ(\zeta)[x,y]$ of the polynomial $g_8$ defined in Corollary \ref{cp1qm}).
\qed

\subsection{Corollary} \label{cp1qm}
The images in $\PP^1_{QM}$ under the parametrization given in 
Proposition \ref{p1qm} of the zeroes of the polynomials 
$$
g_6\,:=\,xy(x^4-y^4),\qquad g_8\,:=\,x^8+14x^4y^4+y^8,\qquad
g_{12}\,:=\,x^{12} - 33x^8y^4 - 33x^4y^8 + y^{12}~,
$$
are the orbits of the points in $\PP^1_{QM}$ with a non-trivial stabilizer in $S_4$. 
Moreover, the rational function
$$
G\,:=\,g_6^4/g_8^3\,:\quad\PP^1_{QM}\,\longrightarrow\,\PP^1\,\cong\,\PP^1_{QM}/S_4
$$
defines the quotient map by $S_4$.

\ts
A non-trivial element $\sigma$ in $S_4\subset Aut(\PP^1_{QM})$ 
has two fixed points, 
corresponding to the eigenlines of any lift of $\sigma$ to $GL(2,\CC)$.
The fixed points of $\sigma^k$ are the same as those of $\sigma$ whenever $\sigma^k$ is not the identity on $\PP^1_{QM}$. One now easily verifies that the fixed points of cycles of 
order $3,4,2$ are the zeroes of $g_6,g_8,g_{12}$ respectively.

The quotient map $\PP^1_{QM}\rightarrow \PP^1_{QM}/S_4\cong\PP^1$
has degree $24$.
The rational function $G:=g_6^4/g_8^3$ is $S_4$-invariant and 
defines a map of degree $24$ from $\PP^1_{QM}$ to $\PP^1$, 
hence the quotient map is  given by $G$. 
\qed

\section{The principal polarization}

\subsection{Introduction} In the previous section we considered Abelian surfaces whose endomorphism ring contains $\cO_6$ endowed with a $(1,2)$-polarization.
Rotger proved that an Abelian surface   
whose endomorphism ring is $\cO_6$  admits a unique principal polarization,
\cite[section 7]{Rotger}.
As such a surface is simple, it is the Jacobian of a genus two curve.
The Abel-Jacobi image of the genus two curve provides the principal polarization.
In this section we find the image of such a curve in the Kummer surface in the $(2,4)$-embedding. This allows us to relate these genus two curves to the ones described by Hashimoto and Murabayashi in \cite{HM} in Section \ref{isHM}.

Moreover, we also find an explicit projective model of a surface in the moduli space
$M_{2,4}$ which parametrizes $(2,4)$-polarized Abelian surfaces 
whose endomorphism ring contains $\ZZ[\sqrt{2}]$, see Section \ref{Humbert}.

\subsection{Polarizations} 
To explain how we found genus two curves in the $(2,4)$-polarized Kummer surfaces parametrized by $\PP^1_{QM}$,
it is convenient to first consider the Jacobian $A=Pic^0(C)$ of one of the genus two curves  
given in \cite[Theorem 1.3]{HM}. 
In \cite[Section 3.1]{HM} one finds an explicit description of the principal polarization $E$ and the maximal order $\cO_6$ of $End(A)\cong B_6\cong (-6,2)_\QQ$. 
The element 
$\eta:=(-1+i)/2+k/4\in \cO_6$ has order three, $\eta^3=1$
(with $i^2=-6,j^2=2,k=ij=-ji$). 
We use the same notation for the endomorphism defined by this element. 
Then $\eta^*E$ is again a principal polarization and we obtain a polarization $E'$ which is invariant under $\eta$ as follows:
$$E':=E+\eta^*E+(\eta^2)^*E,\qquad\mbox{with}\quad 
E(\alpha,\beta):=Tr(\rho_1^{-1}\alpha\beta')
$$  
(here we identify the lattice in $\CC^2$ defining $A$ with $\cO_6$
and $\beta\mapsto \beta'$ is the canonical involution on $B_6$).
An explicit computation shows that $E'=3E''$ and that $E''$ defines a polarization of type $(1,2)$ on $A$ and $\eta^*E''=E''$.

Considering $E$ as a class in $H^2(A,\ZZ)$, one has $E^2=2$, since $E$ is a principal polarization. As $\eta$ is an automorphism of $A$ we also get 
$(\eta^*E)^2=((\eta^2)^*E)^2=2$ and 
$E\cdot(\eta^*E)=(\eta^*E)\cdot((\eta^2)^*)E=((\eta^2)^*)E\cdot E$.
Then one finds that $(E')^2=6+2\cdot 3\cdot E\cdot(\eta^*E)$ and as $E'$ defines a polarization of type $(3,6)$ we have $(E')^2=2\cdot3\cdot6=36$, hence $E\cdot(\eta^*E)=5$. Moreover, one find that
$$
E\cdot E''\,=\,E\cdot(E+\eta^*E+(\eta^2)^*E)/3\,=\,(2+5+5)/3\,=\,4~.
$$

Identify the Jacobian of the genus two curve $C$ with $Pic^0(C)=A$ and
identify $C$ with its image image under the Abel-Jacobi map $C\rightarrow Pic^0(C)$, $p\mapsto p-p_0$, where $p_0$ is a Weierstrass point.  If the hyperelliptic involution interchanges the points $q,q'\in C$, then $q+q'$ and $2p_0$
are linearly equivalent and thus $q-p_0=-(q'-p_0)$. 
Hence the curve $C\subset Pic^0(C)$ is symmetric: $(-1)^*C=C$. 
If $p_1,\ldots,p_5$ are the other Weierstrass points of $C$, 
then $2p_i$ is linearly equivalent to $2p_0$, hence 
the five points $p_i-p_0\in C\subset A$, $i=1,\ldots, 5$ are points of order two in $A$.

Let now $\cL$ be a symmetric line bundle on $A$ defining the 
$(1,2)$-polarization $E''$ on $A$. 
As $E\cdot E''=4$, the restriction of $\cL$ to $C$ has degree $4$ and thus
$\cL^{\otimes 2}$ restricts to a degree $8$ line bundle on $C$.
The map given by the even sections $H^0(A,\cL^{\otimes 2})_+$ defines a $\mbox{2:1}$ map from $A$ onto the Kummer surface $A/\pm1$ of $A$ in $\PP^5$. 
As $(2E'')^2=16$, this Kummer surface has degree $16/2=8$. In fact, Barth shows that the Kummer surface is the complete intersection of three quadrics, see Section \ref{redhs}. 
The symmetry of $C$ implies that this image is a rational curve
and the degree of the image of $C$ is four. 
But a rational curve of degree four in a projective space spans at most a $\PP^4$ (since it is embedded by global sections of $\cO_{\PP^1}(4)$, these can be identified with polynomials of degree at most four in one variable).
Moreover, this $\PP^4$ contains at least six of the nodes (the images of the two-torsion points of $A$) of the Kummer surface which lie on $C$.

It should be noticed that any $(2,4)$-polarized Kummer surface in $\PP^5$
contains subsets of four nodes which span only a $\PP^2$ 
(cf.\ \cite[Lemma 5.3]{GS}), these subsets must be avoided to find $C$.

Conversely, given a rational quartic curve on the Kummer surface which passes through exactly 6 nodes, its inverse image in the Abelian surface will be a 
genus two curve $C$. In fact, the general $A$ is simple, hence there are no non-constant maps from a curve of genus at most one to $A$.
The adjunction formula on $A$ shows that $C^2=2$, hence $C$ defines a principal polarization on $A$.  
Rotger \cite[Section 6]{Rotger} proved that an Abelian surface $A$ with $End(A)=\cO_6$ has a unique principal polarization up to isomorphism.
Thus $C$ must be a member of the family of genus two curves in given in
\cite[Theorem 1.3]{HM}.
We summarize the results in this section in the following proposition.
In Proposition \ref{isoc} we determine the curve from \cite{HM} which is isomorphic to $C=C_x$ on the Abelian surface defined by $x\in\PP^1_{QM}$.

\subsection{Proposition}\label{polp}
Let $A$ be an Abelian surface with $\cO_6\subset \mbox{End}(A)$.
Then $A$ has a (unique up to isomorphism)
principal polarization defined by a genus 2 curve $C\subset A$
which is isomorphic to a curve from the family in \cite[Theorem 1.3]{HM}
(see  Section \ref{isHM}).

There is an automorphism of order three $\eta\in Aut(A)$ such that 
$C+\eta^*C+(\eta^2)^*C=3E''$ defines a polarization of type $(3,6)$. 
Let $\cL$ be a symmetric line bundle with $c_1(\cL)=E''$. Then the image of $C$,
symmetrically embedded in $A$,
under the map $A\rightarrow\PP^5$ defined by the subspace 
$H^0(A,\cL^{\otimes 2})_+$, is a rational curve of degree four which passes through exactly six nodes of the Kummer surface of $A$ 
which lie in a hyperplane in $\PP^5$.

Conversely, the inverse image in $A$ of a rational curve which passes through exactly six nodes of the Kummer surface of $A$ is a genus two curve which defines a principal polarization on $A$.

\subsection{A reducible hyperplane section}\label{redhs}
Now we give a hyperplane $H_x\subset\PP^5$ which cuts the Kummer surface $K_x$ for $x\in\PP^1_{QM}$ in two rational curves of degree four, the curves intersect in six points which are nodes of $K_x$.

A general point $x=(x_1:\ldots:x_6)\in M_{2,4}\subset \PP^5$ 
defines a $(2,4)$-polarized Kummer surface $K_x$
which is the complete intersection of the following three quadrics in $X_1,\ldots,X_6$:
{\renewcommand{\arraystretch}{1.3}
$$
\begin{array}{rcl}
q_1&:=&
(x_1^2 + x_2^2)(X_1^2 +X_2^2) \,-\, (x_3^2 + x_4^2)(X_3^2 +X_4^2) \,
-\,(x_5^2 + x_6^2)(X_5^2 + X_6^2)~,\\
q_2&:=&(x_1^2 - x_2^2)(X_1^2 -X_2^2)\, -\, (x_3^2 - x_4^2)(X_3^2 -X_4)^2\, 
-\,(x_5^2 - x_6^2)(X_5^2 - X_6^2)~,\\
q_3&:=&x_1x_2X_1X_2\,- \,x_3x_4X_3X_4  \,-\,x_5x_6X_5X_6~,
\end{array}
$$
} 
(\cite[Proposition 4.6]{Barth}, we used the formulas from \cite[p.68]{Barth} to replace the $\lambda_i,\mu_i$ by the $x_i$, but notice that the factors `$2$' in the formulas for $\lambda_i\mu_i$ should be omitted, so $\lambda_1\mu_1=x_3^3+x_4^2$ etc.).
The $16$ nodes of the Kummer surface are the orbit of $x$ under the action of $T(2,4)[2]$, that is, it is the set 
$$
\mbox{Nodes}(K_x)\,=\,\{p_{a,b,c,d}\,:=\,(\tilde{\sigma}_1^a\tilde{\sigma_2}^{2b}\tilde{\tau}_1^c\tilde{\tau_2}^{2d})(x),
\quad a,b,c,d\in\{0,1\}\,\}~,
$$
cf.\ Section \ref{hact}.
 We considered the following six nodes:
$$
p_{0,0,0,0},\quad p_{0,0,1,1},\quad  p_{0,1,0,0},\quad
p_{0,1,1,0},\quad p_{1,1,1,0},\quad p_{1,1,1,1}.
$$
For general $x\in P^1_{QM}$ one finds that these six nodes span only a hyperplane $H_x$ in $\PP^5$.

Using Magma we found that over the quadratic extension of the function field 
$\QQ(\zeta)(u)$ of $\PP^1_{QM}$ (where $\zeta^4=-1$ and $u=x/y$) defined by 
$w^2=u^8+14u^4+1$, the intersection of $H_x$ and $K_x$ is reducible and consists of two rational curves of degree four, meeting in the 6 nodes.

We parametrize $H_x$ by $t_1p_{0,0,0,0}+\ldots+t_5p_{1,1,1,0}$. Then
Magma shows that the rational function $t_4/t_5$ 
restricted to each of the two components is a generator of the function field of each of the two components. Thus $t_4/t_5$ provides a coordinate on each component and, for each component, we computed the value (in $\PP^1=\CC\cup\{\infty\}$) 
of the coordinate in the 6 nodes.
The genus two curve $C=C_x$ is the double cover of $\PP^1$ branched 
in these six points.

\subsection{Invariants of genus two curves}
A genus two curve over a field of characteristic $0$ defines a homogeneous sextic polynomial in two variables, uniquely determined up to the action of $Aut(\PP^1)$.
In \cite[p.620]{Igusa}, Igusa defines invariants $A,B,C,D$  of a sextic and defines
further invariants $J_i$, $i=2,4,6,10$, as follows \cite[p.621-622]{Igusa}:
$$
J_2=2^{-3}A,\quad J_4=2^{-5}3^{-1}(4J_2^2-B),\quad 
J_6=2^{-6}3^{-2}(8J_2^3-160J_2J_4-C),\quad
J_{10}=2^{-12}D~.
$$
In \cite[Theorem 6]{Igusa}, Igusa showed that the moduli space of genus two curves over 
$Spec(\ZZ)$ is a (singular) affine scheme which can be embedded in the affine space 
$\bA^{10}_\ZZ$. Its restriction to $Spec(\ZZ[1/2])$ can be embedded into 
$\bA^{8}_{\ZZ[1/2]}$
using the functions (\cite[p.642]{Igusa})
$$
J_2^5J_{10}^{-1},\quad J_2^3J_4J_{10}^{-1},\quad J_2^3J_4^2J_{10}^{-1},\quad
J_2^2J_6J_{10}^{-1},\quad J_4J_6J_{10}^{-1},\quad J_2J_6^3J_{10}^{-2},\quad
J_4^5J_{10}^{-2},\quad J_6^5J_{10}^{-3}~.
$$

From this one finds that
over $Spec(\QQ)$ one can embed the moduli space into $\bA^8_\QQ$ 
using 8 functions $i_1\ldots,i_8$ as above but with 
$J_2,\ldots,J_{10}$ replaced by 
$A,\ldots,D$. 
In  case $A\neq 0$, one can use the three regular functions 
$$
j_1\,:=\,A^5/D,\qquad j_2\,:=\,A^3B/D,\qquad j_3\,:=\,A^2C/D
$$
to express $i_1,\ldots,i_8$ as
$$
j_1,\quad j_2,\quad j_2^2/j_1,\quad j_3,\quad j_2j_4/j_1,\quad j_4^3/j_1,\quad
j_2^5/j_1^3,\quad j_4^5/j_1^2~.
$$
Thus the open subset of the moduli space over $\QQ$ where $A\neq 0$ can be embedded in $\bA^3_\QQ$ using these three functions. In particular, two homogeneous sextic polynomials $f,g$ with complex coefficients and with $A(f),A(g)\neq 0$ define isomorphic genus two curves over $\CC$ if and only if $j_i(f)=j_i(g)$ for $i=1,2,3$ (see also \cite{Mestre},\cite{CQ}).

\subsection{Invariants of the curve $C_x$}\label{invCx}
With the Magma command `IgusaClebschInvariants' we computed the invariants
for each of the two genus curves which are the double covers of the two rational curves in $H_x\cap K_x$, they turn out to be isomorphic as expected from Rotger's uniqueness result.
We denote by $C_x$ the corresponding genus two curve.
For the general $x\in \PP^1_{QM}$ the invariant $A=A(C_x)$ is non-zero and 
$$
j_1(C_x)\,=\,-3^52^{-5}\frac{(1-64G(x))^5}{G(x)^3},\quad
j_2(C_x)\,=\,3^52^{-3}\frac{(1-64G(x))^3}{G(x)^2},
$$
and
$$
j_3(C_x)\,=\,3^42^{-3}\frac{(1-64G(x))^2(1-80G(x))}{G(x)^2}~.
$$
Notice that the invariants are rational functions in the $S_4$-invariant function
$G=g_6^4/g_8^3$ on $\PP^1_{QM}$, as expected. 
Moreover, the $j_i(C_x)$ actually determine $G(x)$:
$$
G(x)\,=\,\frac{(j_2(x)/j_3(x))-3}{80(j_2(x)/j_3(x))-192}~,
$$
hence the classifying map from (an open subset of) $\PP^1_{QM}/S_4$ 
to the moduli space of genus two curves is a birational isomorphism onto its image.

\subsection{The genus two curves from Hashimoto-Murabayashi}\label{isHM}
In \cite[Theorem 1.3]{HM}, Hashimoto and Murabayashi determine an explicit family of genus two curves $C_{s,t}$ whose Jacobians have quaternionic multiplication by 
the maximal order $\cO_6$.
They are parametrized by the elliptic curve
$$
E_{HM}:\qquad g(t,s)\,=\,4s^2t^2\,-\,s^2\,+\,t^2\,+\,2\,=\,0~.
$$
Using the following rational functions on this curve:
$$
P\,:=\,-2(s+t),\quad R\,:=\,-2(s-t),\quad 
Q\,:=\,\frac{(1+2t^2)(11-28t^2+8t^4)}{3(1-t^2)(1-4t^2)}~,
$$
the genus two curve $C_{s,t}$ corresponding to the point $(s,t)\in E_{HM}$ 
is defined by the Weierstrass equation:
$$
C_{s,t}:\qquad Y^2\,=\,X(X^4\,-\,PX^3\,+\,QX^2\,-\,RX\,+\,1)~.
$$
By the unicity result from \cite[section 7]{Rotger} we know that this one parameter family 
of genus two curves should be the same as the one parametrized by $\PP^1_{QM}$.
Indeed one has:

\subsection{Proposition} \label{isoc}
The genus two curve $C_x$ defined by $x\in\PP^1_{QM}$ is isomorphic to the curve 
$C_{s,t}$ if and only if $G(x)=H(t)$ (so the isomorphism class of $C_{s,t}$ does not depend on $s$) where
$$
H(t)\,:=
\,\frac{4(t-1)^2(t+1)^2(t^2+1/2)^4}
{27((1-2t)(1+2t))^3}~.
$$

\ts
This follows from a direct Magma computation of the invariants $j_i$ for the $C_{s,t}$.
In particular, the classifying map of the Hashimoto-Murabayashi family 
has degree $12$ on the $t$-line (and degree $6$ on the $u:=t^2$-line), 
and this degree six cover is not Galois. 
\qed

\subsection{Special points} 
In Section \ref{p1qm} we observed that $S_4$ acts on $\PP^1_{QM}$ and has three orbits which have less then $24$ elements, they are the zeroes of the polynomials $g_d$, of degree $d$, with $d=6,8,12$. 
In case $d=12$ one finds that for example $x=\zeta$ is a zero of $g_{12}$. The invariants $j_i(C_x)$ are the same as the invariants of the curve $C_{s,t}$ from \cite{HM} with
$(t,s)=(0,\sqrt{2})$. In \cite[Example 1.5]{HM} one finds that the Jacobian of this curve is isogenous to a product of two elliptic curves with complex multiplication by $\ZZ[\sqrt{-6}]$. 

In case $d=6,8$ one finds that the invariants $j_i(C_x)$ are infinite,
hence these points do not correspond to Jacobians of genus two curves but to products of
two elliptic curves (with the product polarization). In case $g_6(x)=0$ one finds that the intersection of the plane $H_x$ with the Kummer surface $K_x$ consists of four conics, each of which passes through four nodes (and there are now $8$ nodes in $H_x\cap K_x$).
The inverse image of each conic in the Abelian surface $A_x$ is an elliptic curve which is isomorphic to $E_4:=\CC/\ZZ[i]$, and one finds that $A_x\cong E_4\times E_4$, but the $(1,2)$ polarization is not the product polarization. 
The point $(t,s)=(\sqrt{-2}/2,\sqrt{2}/2)\in E_{HM}$ defines the same point in the Shimura curve
$\PP^1_{QM}/S_4$ as the zeroes of $g_6$. It corresponds to
the degenerate curve $C_{t,s}$ in \cite[Example 1.4]{HM}, which has a normalization which is
isomorphic to $E_4$.

In case $d=8$ one has $A_x\cong E_3\times E_3$ and, with the $(1,2)$-polarization, it is the 
surface $A_3$ that we defined in Section \ref{prelc}.
According to \cite[Theorem 4.9]{Barth} a point 
$x\in M_{2,4}$ defines an Abelian surface $A_x$ if and only if $r(x)\neq 0$ where
$r=r_{12}r_{13}r_{23}$ is defined in \cite[Proposition 3.2]{Barth} (the $r_{jk}$ are polynomials in $\lambda_i^2,\mu_i^2$ and these again can be represented by polynomials in the $x_i$, see \cite[p. 68]{Barth}. One can choose these polynomials as follows:
$$
r_{12}=-4r_{13}=-4r_{23}\,=\,16(x_1x_6 - x_2x_5)(x_1x_6 + x_2x_5)(x_1x_5 - x_2x_6)(x_1x_5 +x_2x_6)~,
$$
and thus $r=16r_{12}^3$.
Restricting $r$ to $\PP^1_{QM}$ and pulling back along the parametrization to $\PP^1$, one finds that 
$
r\,=\,c g_8^3
$,
where $g_8$ is as in Section \ref{p1qm} and $c$ is a non-zero constant.
More in general, we have the following result.

\subsection{Proposition} 
The image of the period matrices $\tau\in\HH_2$ with $\tau_{12}=\tau_{21}=0$
in $M_{2,4}\subset\PP^5$ is the intersection of $M_{2,4}$ with the Segre threefold
which is the image of the map
$$
S_{1,2}:\,
\PP^1\times\PP^2\,\longrightarrow\,\PP^5,\qquad 
\big((u_0:u_1),(w_0:w_1,w_2)\big)\,\longrightarrow\,(x_1:\ldots:x_6)
$$
where the coordinate functions are
$$
\begin{array}{lll}
x_1=u_0w_0,\quad &x_3=u_0w_1,\quad &x_5=u_0w_2~,\\
x_2=u_1w_0,\quad &x_4=u_1w_1,\quad &x_6=u_1w_2~.
\end{array}
$$
The image of $S_{1,2}$ intersects $\PP^1_{QM}$ in two points which are zeroes of $g_8$.
Moreover, the surface $S_{1,2}(\PP^1\times\PP^2)\, \cap \,M_{2,4}$
is an irreducible component of $(r=0)\,\cap\, M_{2,4}$.

\ts If $\tau_{12}=\tau_{21}=0$, then by looking at the Fourier series which define the theta constants, one finds that
$\vartheta[{}^{ab}_{00}](\tau)=
\vartheta[{}^a_0](\tau_{11})\vartheta[{}^b_0](\tau_{22})$. 
The definition of the $x_i$'s in terms of the standard delta functions in $V(2,4)$,
$u_nv_m=\vartheta[{}^{ab}_{00}](\tau)$ with $(a,b)=(n/2,m/4)$ (\cite[p.53]{Barth}), 
then shows that 
the map $\HH_2\rightarrow\PP^5$ restricted to these period matrices 
is the composition of the map
$$
\HH_1\times\HH_1\,\longrightarrow\,\PP^1\times\PP^2,\qquad
(\tau_1,\tau_2)\,\longmapsto\,
$$
$$
\big((\vartheta[{}^0_0](\tau_{1}):\vartheta[{}^b_0](\tau_{1})),
(\vartheta[{}^0_0](\tau_{2})+\vartheta[{}^b_0](\tau_{2}):
\vartheta[{}^a_0](\tau_{2})+\theta[{}^c_0](\tau_{2}):
\theta[{}^0_0](\tau_{2})-\theta[{}^b_0](\tau_{2}))\big)
$$
with the Segre map as above and $a,b,c=1/4,1/2,3/4$ respectively.

The ideal of the image of $S_{1,2}$ is generated by three quadrics, restricting these to
$P^1_{QM}$ one finds that the intersection of the image with $P^1_{QM}$ is defined
by the quadratic polynomial
$x^2 + (\zeta^2 - 1)xy + \zeta^2y^2$, which is a factor of $g_8$.

The factor $x_1x_6 - x_2x_5$ of $r$ is in the ideal of $S_{1,2}(\PP^1\times\PP^2)$,
hence this surface is an irreducible component of $(r=0)\,\cap\, M_{2,4}$.
\qed

\subsection{Remark}
The intersection of the image of $S_{1,2}$ with $M_{2,4}$, which is defined by $f_1=f_2=0$ (cf.\ Section \ref{m24}), is the image of the surface
$$
\PP^1\times C_F,\quad(\subset\PP^1\times\PP^2),\qquad
C_F\,:\quad w_0^4\,-\,w_1^4\,-\,w_2^4\,=\,0~.
$$
The curves $\PP^1$ and $C_F$ here are both elliptic modular curves (defined by the totally symmetric theta structures associated to the divisors $2O$ and $4O$, where $O$ is the origin of the elliptic curve).

\subsection{A Humbert surface}\label{Humbert}
In section \ref{redhs} we considered six nodes of the Kummer surface $K_x$,
$p_{0,0,0,0},p_{0,0,1,1},\ldots,p_{1,1,1,1}$, which had the property that 
for a general $x\in \PP^1_{QM}$ these six nodes span only a hyperplane in $\PP^5$.
For general $x\in M_{2,4}$ however these nodes do span all of $\PP^5$.
They span at most a hyperplane if the determinant $F$ of the $6\times 6$ matrix 
whose rows are the homogeneous coordinates of the nodes, is equal to zero.
$$
F\,=\,\det\,\left(\begin{array}{rrrrrr}
  x_1  &x_2  &x_3  &x_4  &x_5  &x_6\\
-x_2  &x_1  &x_4 &-x_3  &x_6 &-x_5\\
-x_2  &x_1 &-x_4  &x_3  &x_6 &-x_5\\
  x_1 &-x_2  &x_3 &-x_4 &-x_5  &x_6\\
  x_1  &x_2  &x_3  &x_4 &-x_5 &-x_6\\
  x_1 &-x_2 &-x_3  &x_4  &x_5 &-x_6\\
  \end{array}\right)\;=\,16(x_1^2x_3^2x_5x_6\, +\,\ldots\,
-x_2^2x_4^2x_5x_6)~.
$$
Then $F$ is a homogeneous polynomial of degree six in the coordinates of $x$ which has $8$ terms.
Let $D_F$ be the divisor in $M_{2,4}$ defined by $F=0$, then $\PP^1_{QM}$
is contained in (the support of) $D_F$. 
Magma shows that $D_F$ has $12$ irreducible components,
the only one of these which contains $\PP^1_{QM}$ is the surface $S_2\subset\PP^5$ defined by
$$
S_2:\quad x_1^2 - x_2^2 - x_5^2 - x_6^2\,=\,
x_1x_2 - x_4^2 - x_5x_6\,=\,x_3^2 - x_4^2 - 2x_5x_6\,=\,0~.
$$
Magma verified that $S_2$ is a smooth surface, hence it is a K3 surface. 

\subsection{Proposition} The surface $S_2\subset M_{2,4}$ parametrizes 
Abelian surfaces $A$ with $\ZZ[\sqrt{2}]\subset End(A)$.

\ts
For a general point $x$ in $S_2$, the hyperplane spanned by the six nodes 
intersects $K_x$ in a one-dimensional subscheme which is the complete intersection 
of three quadrics and which has six nodes. The arithmetic genus of a smooth complete intersection of three quadrics in $H_x=\PP^4$ is only five, hence this subscheme must be reducible.
In the case $x\in \PP^1_{QM}$, this subscheme is the union of two
smooth rational curves of degree four intersecting transversally in the six nodes.
Thus for general $x\in S_2$, the intersection must consist also consist of two such rational curves. Let $C\subset A_x$ be the genus two curve in the Abelian surface $A_x$ defined by $x$ which is the inverse image of one of these components. 
Then $C^2=2$ and $C\cdot \cL=4$, where $\cL$ defines the $(1,2)$-polarization. 
Now we apply \cite[Proposition 5.2.3]{BL} 
to the endomorphism $f=\phi_C^{-1}\phi_\cL$ of $A_x$ defined by these polarizations.
We find that the characteristic polynomial of $f$ is $t^2-4t+2$. As its  roots are  $2\pm\sqrt{2}$, we conclude that $\ZZ[\sqrt{2}]\subset End(A_x)$. 
\qed

\

\end{document}